\documentclass{article}%
\usepackage{amsmath}%
\usepackage{amsfonts}%
\usepackage{amssymb}%
\usepackage{graphicx}
\newtheorem{theorem}{Theorem} [section]

\newtheorem{conjecture}[theorem]{Conjecture}
\newtheorem{corollary}[theorem]{Corollary}

\newtheorem{lemma}[theorem]{Lemma}

\newtheorem{proposition}[theorem]{Proposition}

\newenvironment{proof}[1][Proof]{\noindent\textbf{#1.} }{\ \rule{0.5em}{0.5em}}
\begin{document}

\title{Very well-covered graphs and the unimodality conjecture}
\author{Vadim E. Levit and Eugen Mandrescu\\Department of Computer Science\\Holon Academic Institute of Technology\\52 Golomb Str., P.O. Box 305\\Holon 58102, ISRAEL}
\maketitle

\begin{abstract}
If $s_{k}$ denotes the number of stable sets of cardinality $k$ in the graph
$G$, then $I(G;x)=\sum\limits_{k=0}^{\alpha(G)}s_{k}x^{k}$ is the
\textit{independence polynomial} of $G$ (Gutman, Harary, 1983), where
$\alpha(G)$ is the size of a maximum stable set in $G$. Alavi, Malde, Schwenk
and Erd\"{o}s (1987) conjectured that $I(T,x)$ is unimodal for any tree $T$,
while, in general, they proved that for any permutation $\pi$ of
$\{1,2,...,\alpha\}$ there is a graph $G$ with $\alpha(G)=\alpha$ such that
$s_{\pi(1)}<s_{\pi(2)}<...<s_{\pi(\alpha)}$. Brown, Dilcher and Nowakowski
(2000) conjectured that $I(G;x)$ is unimodal for any well-covered graph.
Michael and Traves (2002) provided examples of well-covered graphs with
non-unimodal independence polynomials. They proposed the so-called
\textquotedblright roller-coaster\textquotedblright\ conjecture: for a
well-covered graph, the subsequence $(s_{\left\lceil \alpha/2\right\rceil
},s_{\left\lceil \alpha/2\right\rceil +1},...,s_{\alpha})$ is unconstrained in
the sense of Alavi \textit{et al}. The conjecture of Brown \textit{et al.} is
still open for very well-covered graphs, and it is worth mentioning that,
apart from $K_{1}$ and the chordless cycle $C_{7}$, connected well-covered
graphs of girth $\geq6$ are very well-covered (Finbow, Hartnell and
Nowakowski, 1993).

In this paper we prove that $s_{\left\lceil (2\alpha-1)/3\right\rceil }%
\geq...\geq s_{\alpha-1}\geq s_{\alpha}$ are valid for any \textit{(a)}
bipartite graph $G$ with $\alpha(G)=\alpha$; \textit{(b)} quasi-regularizable
graph $G$ on $2\alpha(G)=2\alpha$ vertices. In particular, we infer that this
is true for \textit{(a)} trees, thus doing a step in an attempt to prove Alavi
\textit{et al.}' conjecture; \textit{(b)} very well-covered graphs.
Consequently, for this case, the unconstrained subsequence appearing in the
roller-coaster conjecture can be shorten to $(s_{\left\lceil \alpha
/2\right\rceil },s_{\left\lceil \alpha/2\right\rceil +1},...,s_{\left\lceil
(2\alpha-1)/3\right\rceil })$. We also show that the independence polynomial
of a very well-covered graph $G$ is unimodal for $\alpha(G)\leq9$, and is
log-concave whenever $\alpha(G)\leq5$.

\textbf{key words:}\textit{\ stable set, independence polynomial, unimodal
sequence, quasi-regularizable graph, bipartite graph, very well-covered
graph.}

\end{abstract}

\section{Introduction}

Throughout this paper $G=(V,E)$ is a simple (i.e., a finite, undirected,
loopless and without multiple edges) graph with vertex set $V=V(G)$ and edge
set $E=E(G).$ If $X\subset V$, then $G[X]$ is the subgraph of $G$ spanned by
$X$. By $G-W$ we mean the subgraph $G[V-W]$, if $W\subset V(G)$. We also
denote by $G-F$ the partial subgraph of $G$ obtained by deleting the edges of
$F$, for $F\subset E(G)$, and we write shortly $G-e$, whenever $F$ $=\{e\} $.

A vertex $v$ is \textit{pendant} if its neighborhood $N(v)=\{u:u\in V,uv\in
E\}$ contains only one vertex; an edge $e=uv$ is \textit{pendant} if one of
its endpoints is a pendant vertex. $\overline{G}$ stands for the complement of
$G$, while $K_{n},P_{n},C_{n}$ denote respectively, the complete graph on
$n\geq1$ vertices, the chordless path on $n\geq1$ vertices, and the chordless
cycle on $n\geq3$ vertices. As usual, a \textit{tree} is an acyclic connected graph.

A set of pairwise non-adjacent vertices is called \textit{stable}. If $S$ is a
stable set, then we denote $N(S)=\{v:N(v)\cap S\neq\emptyset\}$ and
$N[S]=N(S)\cup S$. A stable set of maximum size will be referred to as a
\textit{maximum stable set} of $G$. The \textit{stability number }of $G$,
denoted by $\alpha(G)$, is the cardinality of a maximum stable set in $G$, and
$\omega(G)=\alpha(\overline{G})$.

The \textit{disjoint union} of the graphs $G_{1},G_{2}$ is the graph
$G=G_{1}\sqcup G_{2}$ having as vertex set and edge set the disjoint unions of
$V(G_{1}),V(G_{2})$ and $E(G_{1}),E(G_{2})$, respectively.

If $G_{1},G_{2}$ are disjoint graphs, then their \textit{Zykov sum}, (Zykov,
\cite{Zykov}, \cite{Zykov1}), is the graph $G_{1}+G_{2}$ with
\begin{align*}
V(G_{1}+G_{2})  & =V(G_{1})\cup V(G_{2}),\\
E(G_{1}+G_{2})  & =E(G_{1})\cup E(G_{2})\cup\{v_{1}v_{2}:v_{1}\in
V(G_{1}),v_{2}\in V(G_{2})\}.
\end{align*}
In particular, $\sqcup nG$ and $+nG$ denote the disjoint union and Zykov sum,
respectively, of $n>1$ copies of the graph $G$.

A graph $G$ is called \textit{quasi-regularizable} if one can replace each
edge of $G$ with a non-negative integer number of parallel copies, so as to
obtain a regular multigraph of degree $\neq0$ (see \cite{Berge1},
\cite{Berge2}). Evidently, a disconnected quasi-regularizable graph has no
isolated vertices. Moreover, a disconnected graph is quasi-regularizable if
and only if any of its connected components spans a quasi-regularizable graph.
The following characterization of quasi-regularizable graphs, due to Berge, we
shall use in the sequel.

\begin{theorem}
\cite{Berge1}\label{th1} A graph $G$ is quasi-regularizable if and only if
$\left|  S\right|  \leq\left|  N(S)\right|  $ holds for any stable set $S$ of
$G$.
\end{theorem}

Let $s_{k}$ be the number of stable sets in $G$ of cardinality $k\in
\{0,1,...,\alpha(G)\}$. The polynomial
\[
I(G;x)=s_{0}+s_{1}x+s_{2}x^{2}+...+s_{\alpha}x^{\alpha},\;\alpha=\alpha(G),
\]
is called the \textit{independence polynomial} of $G$ (Gutman and Harary,
\cite{GutHar}). Various properties of this polynomial are presented in a
number of papers, like \cite{GutHar}, \cite{BrownDilNow}, \cite{BrownHickNow},
\cite{HodeLi}, \cite{LevitMan2}, \cite{LevitMan3}, \cite{LevitMan4},
\cite{LevitMan5}, \cite{LevitMan6}, \cite{MichaelTraves}.

A finite sequence of real numbers $(a_{0},a_{1},a_{2},...,a_{n})$ is said to be:

\begin{itemize}
\item \textit{unimodal} if there is some $k\in\{0,1,...,n\}$, called the
\textit{mode} of the sequence, such that $a_{0}\leq...\leq a_{k-1}\leq
a_{k}\geq a_{k+1}\geq...\geq a_{n}$;

\item \textit{log-concave} if $a_{i}^{2}\geq a_{i-1}\cdot a_{i+1}$ holds for
$i\in\{1,2,...,n-1\}$.
\end{itemize}

It is known that any log-concave sequence of positive numbers is also unimodal.

A polynomial is called \textit{unimodal (log-concave)} if the sequence of its
coefficients is unimodal (log-concave, respectively). For instance, the
independence polynomial $I(K_{1,3};x)=1+4x+3x^{2}+x^{3}$ is unimodal. However,
the independence polynomial of $G=K_{24}+(K_{3}\sqcup K_{3}\sqcup K_{4})$ is
not unimodal, since $I(G;x)=1+34x+33x^{2}+36x^{3}$ (for other examples, see
\cite{AlMalSchErdos}). Moreover, Alavi \textit{et al.} \cite{AlMalSchErdos}
proved that for any permutation $\pi$ of $\{1,2,...,\alpha\}$ there is a graph
$G$ with $\alpha(G)=\alpha$ such that $s_{\pi(1)}<s_{\pi(2)}<...<s_{\pi
(\alpha)}$. Nevertheless, for trees, they stated the following still open conjecture.

\begin{conjecture}
\cite{AlMalSchErdos}\label{conj1} The independence polynomial of a tree is unimodal.
\end{conjecture}

A graph $G$ is called \textit{well-covered} if all its maximal stable sets
have the same cardinality, (Plummer, \cite{Plum}). If, in addition, $G$ has no
isolated vertices and its order equals $2\alpha(G)$, then $G$ is \textit{very
well-covered} (Favaron, \cite{Favaron}).

By $G^{*}$ we mean the graph obtained from $G$ by appending a single pendant
edge to each vertex of $G$, \cite{DuttonChanBrigham}. Let us notice that
$G^{*}$ is well-covered (see, for instance, \cite{LevMan0}), and $\alpha
(G^{*})=n$. In fact, $G^{*}$ is very well-covered. Moreover, the following
result shows that, under certain conditions, any well-covered graph equals
$G^{*}$ for some graph $G$.

\begin{theorem}
\cite{FinHarNow}\label{th3} Let $G$ be a connected graph of girth $\geq6$,
which is isomorphic to neither $C_{7}$ nor $K_{1}$. Then $G$ is well-covered
if and only if its pendant edges form a perfect matching.
\end{theorem}

In other words, Theorem \ref{th3} shows that, apart from $K_{1}$ and $C_{7}$,
connected well-covered graphs of girth $\geq6$ are very well-covered. For
example, a tree $T\neq K_{1}$ could be only very well-covered, and this is the
case if and only if $T=G^{\ast}$ for some tree $G$ (see also \cite{Ravindra},
\cite{Favaron}, \cite{LevitMan1}).

In \cite{BrownDilNow} it was conjectured that the independence polynomial of
any well-covered graph $G$ is unimodal. Michael and Traves
\cite{MichaelTraves} proved that this conjecture is true for $\alpha
(G)\in\{1,2,3\}$, but it is false for $\alpha(G)\in\{4,5,6,7\}$. A family of
well-covered graphs with non-unimodal independence polynomials and stability
numbers $\geq8$ is presented in \cite{LevitMan6}. However, the conjecture is
still open for very well-covered graphs. In \cite{LevitMan2} and
\cite{LevitMan3}, unimodality of independence polynomials of a number of
well-covered trees (e.g., $P_{n}^{\ast},K_{1,n}^{\ast}$) was validated, using
the fact that the independence polynomial of a claw-free graph is unimodal (a
result due to Hamidoune, \cite{Hamidoune}). We also showed that $I(G^{\ast
};x)$ is unimodal for any $G^{\ast}$ whose skeleton $G$ has $\alpha(G)\leq4$
(see \cite{LevitMan4}).

Michael and Traves formulated (and verified for well-covered graphs with
stability numbers $\leq7$) the following so-called "\textit{roller-coaster}"\ conjecture.

\begin{conjecture}
\cite{MichaelTraves} For any permutation $\pi$ of the set $\{\left\lceil
\alpha/2\right\rceil ,\left\lceil \alpha/2\right\rceil +1,...,\alpha\}$, there
exists a well-covered graph $G$, with $\alpha(G)=\alpha$, whose sequence
$(s_{0},s_{1},...,s_{\alpha})$ satisfies $s_{\pi\left(  \left\lceil
\alpha/2\right\rceil \right)  }<s_{\pi\left(  \left\lceil \alpha/2\right\rceil
+1\right)  }<...<s_{\pi\left(  \alpha\right)  }$.
\end{conjecture}

In this paper we prove that if $G$ is a quasi-regularizable graph on
$2\alpha(G)$ vertices, then $s_{\left\lceil (2\alpha(G)-1)/3\right\rceil }\geq
s_{\left\lceil (2\alpha(G)-1)/3\right\rceil +1}\geq...\geq s_{\alpha(G)}$,
while if $G$ is a perfect graph, then $s_{\left\lceil (\omega\alpha
-1)/(\omega+1)\right\rceil }\geq s_{\left\lceil (\omega\alpha-1)/(\omega
+1)\right\rceil +1}\geq...\geq s_{\alpha}$, where $\alpha=\alpha
(G),\omega=\omega(G)$. We infer that for very well-covered graphs, the domain
of the roller-coaster conjecture can be shorten to $\{\left\lceil
\alpha/2\right\rceil ,\left\lceil \alpha/2\right\rceil +1,...,\left\lceil
(2\alpha-1)/3\right\rceil \}$. Moreover, we show that the independence
polynomial of a very well-covered graph $G$ is unimodal for $\alpha(G)\leq9 $,
and log-concave, whenever $\alpha(G)\leq5$.

\section{Results}

In \cite{BrownDilNow} it was shown that $s_{k-1}\leq k\cdot s_{k}$ and
$s_{k}\leq(n-k+1)\cdot s_{k-1},1\leq k\leq\alpha(G)$, are true for any
well-covered graph $G$ on $n$ vertices.

\begin{proposition}
\label{prop2}\cite{MichaelTraves}, \cite{LevitMan5} If $G$ is a well-covered
graph with $\alpha(G)=\alpha$, then the following statements are true:

\emph{(i)} $(\alpha-k)\cdot s_{k}\leq(k+1)\cdot s_{k+1}$ holds for $0\leq
k<\alpha$;

\emph{(ii)} $s_{k-1}\leq s_{k}$ for any $1\leq k\leq(\alpha+1)/2$.
\end{proposition}

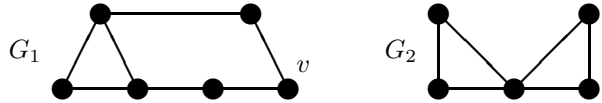
\begin{figure}[h]
\setlength{\unitlength}{1cm}\begin{picture}(5,1.3)\thicklines
\multiput(3,0)(1,0){4}{\circle*{0.29}}
\multiput(3.5,1)(2,0){2}{\circle*{0.29}}
\put(3,0){\line(1,0){3}}
\put(3,0){\line(1,2){0.5}}
\put(3.5,1){\line(1,0){2}}
\put(4,0){\line(-1,2){0.5}}
\put(6,0){\line(-1,2){0.5}}
\put(6.2,0.3){\makebox(0,0){$v$}}
\put(2.5,0.5){\makebox(0,0){$G_{1}$}}
\multiput(8,0)(1,0){3}{\circle*{0.29}}
\multiput(8,1)(2,0){2}{\circle*{0.29}}
\put(8,0){\line(1,0){2}}
\put(8,0){\line(0,1){1}}
\put(8,1){\line(1,-1){1}}
\put(9,0){\line(1,1){1}}
\put(10,0){\line(0,1){1}}
\put(7.5,0.5){\makebox(0,0){$G_{2}$}}
\end{picture}
\caption{Non-well-covered graph.}%
\label{fig3}%
\end{figure}

Notice that Proposition \ref{prop2}\emph{(i)} can fail for non-well-covered
graphs, e.g., the graph $G_{1}$ in Figure \ref{fig3} has $\alpha(G_{1})=3$ and
$(\alpha(G_{1})-2)\cdot s_{2}=8>3=(2+1)\cdot s_{3}$. However, there are
non-well-covered graphs satisfying Proposition \ref{prop2}\emph{(i)}, for
instance, the graph $G_{2}$ in Figure \ref{fig3}. Since $I(G_{1}%
;x)=1+6x+8x^{2}+x^{3}$ and $I(G_{2};x)=1+5x+4x^{2}$, we see that both $G_{1}$
and $G_{2}$ satisfy Proposition \ref{prop2}\emph{(ii)}. On the other hand,
$K_{1,3}$ does not agree with Proposition \ref{prop2}\emph{(ii)}, because
$\alpha(K_{1,3})=3,I(K_{1,3};x)=1+4x+3x^{2}+x^{3}$, while $s_{1}=4>3=s_{2}$.

\begin{corollary}
\label{cor1}If $G$ is a well-covered graph with $\alpha(G)=\alpha$,

then $s_{k}\leq s_{\alpha-k}$ for $0\leq k\leq\alpha/2$.
\end{corollary}

\begin{proof}
Let $k\in\{0,...,\left\lfloor \alpha/2\right\rfloor \}$. According to
Proposition \ref{prop2}\emph{(i)}, we obtain successively, that:
\begin{align*}
(\alpha-k)\cdot s_{k}  & \leq(k+1)\cdot s_{k+1},\\
(\alpha-k-1)\cdot s_{k+1}  & \leq(k+2)\cdot s_{k+2},\\
& ...\\
(k+1)\cdot s_{\alpha-k-1}  & \leq(\alpha-k)\cdot s_{\alpha-k}.
\end{align*}
By multiplying these inequalities, we get
\begin{align*}
& (\alpha-k)(\alpha-k-1)(\alpha-k-2)\cdot...\cdot(k+1)\cdot s_{k}\cdot
s_{k+1}\cdot s_{k+2}\cdot...\cdot s_{\alpha-k-1}\\
& \leq(k+1)(k+2)(k+3)\cdot...\cdot(\alpha-k)\cdot s_{k+1}\cdot s_{k+2}\cdot
s_{k+3}\cdot...\cdot s_{\alpha-k-1}\cdot s_{\alpha-k},
\end{align*}
which clearly leads to $s_{k}\leq s_{\alpha-k}$.
\end{proof}

The above Corollary \ref{cor1} fails for some non-well-covered graphs, e.g.,
$K_{1,3\text{ }}$ has $\alpha(K_{1,3})=3$, while $s_{1}=4>3=s_{2}=s_{\alpha
-1}$. Nevertheless, $P_{5}$ is not a well-covered graph, but $I(P_{5}%
;x)=1+5x+6x^{2}+x^{3}$ shows that $P_{5}$ satisfies Corollary \ref{cor1}.

\begin{figure}[h]
\setlength{\unitlength}{1cm}\begin{picture}(5,1.3)\thicklines
\multiput(5,0)(1,0){4}{\circle*{0.29}}
\multiput(6,1)(1,0){2}{\circle*{0.29}}
\put(5,0){\line(1,0){3}}
\put(6,1){\line(1,0){1}}
\put(6,0){\line(1,1){1}}
\put(7,0){\line(0,1){1}}
\end{picture}
\caption{The graph $K_{3}^{\ast}$.}%
\label{fig1}%
\end{figure}
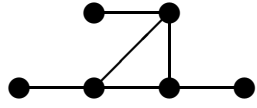

For a graph $G$ of order $n$ and having $\alpha(G)=\alpha$, we denote
\[
\omega_{\alpha-k}=\max\{n-\left\vert N[S]\right\vert :S\ \text{\textit{is a
stable set with} }\left\vert S\right\vert =k\},0\leq k\leq\alpha.
\]
Clearly, $\omega_{0}=0,\omega_{\alpha}=n$. While $\omega_{1}\left(  G\right)
\leq\omega(G)$, it is not necessary that $\omega_{1}\left(  G\right)
=\omega(G)$. For instance, the graph $K_{3}^{\ast}$ (depicted in Figure
\ref{fig1}) has $\omega_{1}=2,\omega(K_{3}^{\ast})=3$. It is worth mentioning
that for any odd chordless cycle $C_{2n+1},n\geq2$, or even chordless path
$P_{2n},n\geq2$, these two parameters are identical. 

\begin{lemma}
\label{lem1}If $G$ is a graph of order $n\geq1$ with $\alpha(G)=\alpha$, then
\[
(k+1)\cdot s_{k+1}\leq\omega_{\alpha-k}\cdot s_{k},0\leq k<\alpha,
\]

in particular, $\alpha\cdot s_{\alpha}\leq\omega_{1}\cdot s_{\alpha-1}%
\leq\omega(G)\cdot s_{\alpha-1}$.
\end{lemma}

\begin{proof}
Let $H=(\mathcal{A},\mathcal{B},\mathcal{W})$ be the bipartite graph defined
as follows: $X\in\mathcal{A}\Leftrightarrow X$ is a stable set in $G$ of size
$k$, then $Y\in\mathcal{B}\Leftrightarrow Y$ is a stable set in $G$ of size
$k+1$, and $XY\in\mathcal{W}\Leftrightarrow X\subset Y$ in $G$. Since any
$Y\in\mathcal{B}$ has exactly $k+1$ subsets of size $k$, it follows that
$\left\vert \mathcal{W}\right\vert =(k+1)\cdot s_{k+1}$. On the other hand, if
$X\in\mathcal{A}$ and, then $X\cup\{v\}\in\mathcal{B}$ for any $v\in
V(G)-N[X]$, i.e., $X$ has at most $\omega_{\alpha-k}$ neighbors in
$\mathcal{B}$. Hence, we get that $(k+1)\cdot s_{k+1}=\left\vert
\mathcal{W}\right\vert \leq\omega_{\alpha-k}\cdot\left\vert \mathcal{A}%
\right\vert =\omega_{\alpha-k}\cdot s_{k}$. In particular, for $k=\alpha-1$,
we obtain $\alpha\cdot s_{\alpha}\leq\omega_{1}\cdot s_{\alpha-1}\leq
\omega(G)\cdot s_{\alpha-1}$.
\end{proof}

Let us remark that there are quasi-regularizable graphs with non-unimodal
independence polynomials, e.g.,

\emph{(a)} $G=K_{10}+\sqcup6K_{1}$ is connected and has
\[
I(G;x)=\left(  1+x\right)  ^{6}+10x=1+\mathbf{16}x+15x^{2}+\mathbf{20}%
x^{3}+15x^{4}+6x^{5}+x^{6};
\]

\emph{(b)} $G=\left(  K_{24}+\sqcup6K_{1}\right)  \sqcup\left(  K_{25}%
+\sqcup6K_{1}\right)  $ is disconnected and has
\begin{align*}
I(G;x)  & =\left(  \left(  1+x\right)  ^{6}+24x\right)  \left(  \left(
1+x\right)  ^{6}+25x\right) \\
& =1+61x+\mathbf{960}x^{2}+955x^{3}+\mathbf{1475}x^{4}+1527x^{5}\\
& +1218x^{6}+841x^{7}+495x^{8}+220x^{9}+66x^{10}+12x^{11}+x^{12}.
\end{align*}

\begin{proposition}
\label{prop1}If $G$ is a quasi-regularizable graph of order $n=2\alpha
(G)=2\alpha$, then

\emph{(i) }$\omega_{\alpha-k}\leq2(\alpha-k),0\leq k\leq\alpha$;

\emph{(ii)} $(k+1)\cdot s_{k+1}\leq2(\alpha-k)\cdot s_{k},0\leq k<\alpha$;

\emph{(iii)} $s_{\left\lceil (2\alpha-1)/3\right\rceil }\geq...\geq
s_{\alpha-1}\geq s_{\alpha}$.
\end{proposition}

\begin{proof}
\emph{(i) }Let $S$ be a stable set in $G$ of size $k\geq0$. According to
Theorem \ref{th1}, it follows that $\left\vert S\right\vert \leq\left\vert
N(S)\right\vert $, which implies $2\cdot\left\vert S\right\vert \leq\left\vert
S\cup N(S)\right\vert =\left\vert N[S]\right\vert $ and, hence, $2\cdot
(\alpha-k)=2\cdot(\alpha-\left\vert S\right\vert )\geq n-\left\vert
N[S]\right\vert $, because $n=2\alpha$. Consequently, we infer that
$\omega_{\alpha-k}\leq2(\alpha-k)$.

\emph{(ii) }The result follows by\emph{\ }combining Lemma \ref{lem1} and part
\emph{(i)}.

\emph{(iii) }The fact that $(k+1)\cdot s_{k+1}\leq2(\alpha-k)\cdot s_{k}$
implies that $s_{k+1}\leq s_{k}$ holds for $k+1\geq2(\alpha-k)$, i.e., for
$k\geq(2\alpha-1)/3$.
\end{proof}

There are no quasi-regularizable graphs $G$ of order $n>2\alpha(G)$ that
satisfy Proposition \ref{prop1}\emph{(i),(ii)}, since for $k=0$, each of them
demands $n\leq2\alpha\left(  G\right)  $.

In addition, for the graphs $G_{1},G_{2}$ in Figure \ref{fig22},
$I(G_{1};x)=1+6x+8x^{2}$ and $I(G_{2};x)=\allowbreak1+8x+19x^{2}+12x^{3}$ show
that Proposition \ref{prop1}\emph{(iii)} is sometimes, but not always, valid
for a quasi-regularizable graph $G$ on $n>2\alpha(G)$ vertices. Notice that
$G_{1}$ is also well-covered, but not very well-covered. 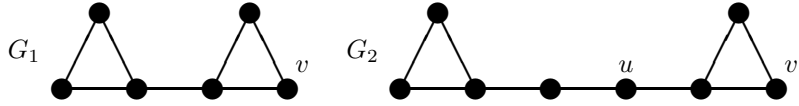
\begin{figure}[h]
\setlength{\unitlength}{1cm}\begin{picture}(5,1)\thicklines
\multiput(1.5,0)(1,0){4}{\circle*{0.29}}
\multiput(2,1)(2,0){2}{\circle*{0.29}}
\put(1.5,0){\line(1,0){3}}
\put(1.5,0){\line(1,2){0.5}}
\put(3.5,0){\line(1,2){0.5}}
\put(2.5,0){\line(-1,2){0.5}}
\put(4.5,0){\line(-1,2){0.5}}
\put(4.7,0.3){\makebox(0,0){$v$}}
\put(1,0.5){\makebox(0,0){$G_{1}$}}
\multiput(6,0)(1,0){6}{\circle*{0.29}}
\multiput(6.5,1)(4,0){2}{\circle*{0.29}}
\put(6,0){\line(1,0){5}}
\put(6,0){\line(1,2){0.5}}
\put(10,0){\line(1,2){0.5}}
\put(7,0){\line(-1,2){0.5}}
\put(11,0){\line(-1,2){0.5}}
\put(11.2,0.3){\makebox(0,0){$v$}}
\put(9,0.3){\makebox(0,0){$u$}}
\put(5.5,0.5){\makebox(0,0){$G_{2}$}}
\end{picture}
\caption{$G_{1},G_{2}$ are quasi-regularizable graphs, but only $G_{1}$ is
well-covered.}%
\label{fig22}%
\end{figure}

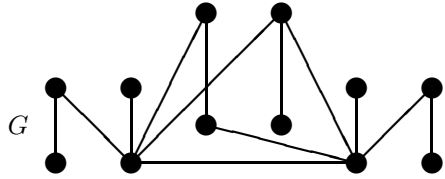
\begin{figure}[h]
\setlength{\unitlength}{1cm}\begin{picture}(5,2)\thicklines
\multiput(4,0)(1,0){2}{\circle*{0.29}}
\multiput(4,1)(1,0){2}{\circle*{0.29}}
\multiput(8,0)(1,0){2}{\circle*{0.29}}
\multiput(8,1)(1,0){2}{\circle*{0.29}}
\multiput(6,2)(1,0){2}{\circle*{0.29}}
\multiput(6,0.5)(1,0){2}{\circle*{0.29}}
\put(4,0){\line(0,1){1}}
\put(4,1){\line(1,-1){1}}
\put(5,0){\line(0,1){1}}
\put(5,0){\line(1,0){3}}
\put(5,0){\line(1,2){1}}
\put(5,0){\line(1,1){2}}
\put(6,2){\line(0,-1){1.5}}
\put(6,0.5){\line(4,-1){2}}
\put(7,2){\line(0,-1){1.5}}
\put(7,2){\line(1,-2){1}}
\put(8,0){\line(0,1){1}}
\put(8,0){\line(1,1){1}}
\put(9,0){\line(0,1){1}}
\put(3.5,0.5){\makebox(0,0){$G$}}
\end{picture}
\caption{A very well-covered graph with a log-concave independence
polynomial.}%
\label{fig60}%
\end{figure}

The graph $G$ in Figure \ref{fig60} is very well-covered and its independence
polynomial $I(G;x)=1+12x+52x^{2}+110x^{3}+123x^{4}+70x^{5}+16x^{6}$ is not
only unimodal but log-concave, as well.

\begin{theorem}
If $G$ is a very well-covered graph of order $n\geq2$ with $\alpha(G)=\alpha$, then

\emph{(i)} $(\alpha-k)\cdot s_{k}\leq(k+1)\cdot s_{k+1}\leq2(\alpha-k)\cdot
s_{k},0\leq k<\alpha$;

\emph{(ii)} $s_{0}\leq s_{1}\leq...\leq s_{\left\lceil \alpha/2\right\rceil }$
and $s_{\left\lceil (2\alpha-1)/3\right\rceil }\geq...\geq s_{\alpha-1}\geq
s_{\alpha}$;

\emph{(iii) }$s_{\alpha-2}\cdot s_{\alpha}\leq s_{\alpha-1}^{2}$, where
$\alpha\geq2$;

\emph{(iv) }$I(G;x)$ is unimodal, while $\alpha\leq9$;

\emph{(v) }$I(G;x)$ is log-concave, while $\alpha\leq5$.
\end{theorem}

\begin{proof}
\emph{(i)} It follows from Proposition \ref{prop2}\emph{(i)} and Proposition
\ref{prop1}, because any well-covered graph without isolated vertices is
quasi-regularizable (see Berge, \cite{Berge1}, \cite{Berge2}).

\emph{(ii) }The assertion is established according to Proposition
\ref{prop2}\emph{(ii)} and Proposition \ref{prop1}.

\emph{(iii) }Taking $k=\alpha-2$ in Proposition \ref{prop2}\emph{(i)}, we get
$2\cdot s_{\alpha-2}\leq(\alpha-1)\cdot s_{\alpha-1}$, while substituting
$k=\alpha-1$ in part \emph{(i)} assures that $\alpha\cdot s_{\alpha}\leq2\cdot
s_{\alpha-1}$, which together lead to $2\alpha\cdot s_{\alpha-2}\cdot
s_{\alpha}\leq2(\alpha-1)\cdot s_{\alpha-1}^{2}$ and, hence, $s_{\alpha
-2}\cdot s_{\alpha}\leq s_{\alpha-1}^{2}$.

\emph{(iv) }By part \emph{(ii)}, $s_{0}\leq s_{1}\leq...\leq s_{\left\lceil
\alpha/2\right\rceil }$ and $s_{\left\lceil (2\alpha-1)/3\right\rceil }%
\geq...\geq s_{\alpha-1}\geq s_{\alpha}$. In addition, the fact that
$\alpha(G)\leq9$ ensures that either$\left\vert \left\lceil \alpha
/2\right\rceil -\left\lceil \left(  2\alpha-1\right)  /3\right\rceil
\right\vert \leq1$.

\emph{(v) }Notice that $s_{0}\cdot s_{2}=\left\vert E(\overline{G})\right\vert
\leq\left\vert V(G)\right\vert ^{2}=s_{1}^{2}$ is true for any graph $G$ with
$\alpha(G)=\alpha\geq2$. By part \emph{(iii)}, $s_{\alpha-2}\cdot s_{\alpha
}\leq s_{\alpha-1}^{2}$. Therefore, we have to check that $s_{k-1}\cdot
s_{k+1}\leq s_{k}^{2}$ only for $k\in\{2,\alpha-2\}$.

Part \emph{(i)} implies that $(\alpha-k+1)\cdot s_{k-1}\leq k\cdot s_{k}$ and
$(k+1)\cdot s_{k+1}\leq2(\alpha-k)\cdot s_{k}$, which together give
\[
(\alpha-k+1)\cdot(k+1)\cdot s_{k-1}s_{k+1}\leq2(\alpha-k)\cdot k\cdot
s_{k}^{2}.
\]
If $(\alpha-k+1)\cdot(k+1)\geq2(\alpha-k)\cdot k$, then $s_{k-1}\cdot
s_{k+1}\leq s_{k}^{2}$. In other words, we are interested to know when
$k^{2}-\alpha k+\alpha+1\geq0$, while $2\leq k\leq\alpha-2$. Since the roots
of $k^{2}-\alpha k+\alpha+1$ are $k_{1,2}=(\alpha\pm\sqrt{\alpha^{2}%
-4\alpha-4})/2$, we conclude the following, depending on $\alpha$:

\emph{(a)} $\alpha\leq4$, then $\alpha^{2}-4\alpha-4<0$ and, hence,
$k^{2}-\alpha k+\alpha+1\geq0$ is valid for any $k$;

\emph{(b)} $\alpha=5$, then $k_{1}=2,k_{2}=3$, and $k^{2}-\alpha
k+\alpha+1\geq0$ is still true for any $k$;

\emph{(c)} $\alpha\geq6$, then $k^{2}-\alpha k+\alpha+1\geq0$ only for $k=1$
and $k=\alpha-1$, because $2<(\alpha-\sqrt{\alpha^{2}-4\alpha-4})/2<4$ and
$2(\alpha-2)<(\alpha+\sqrt{\alpha^{2}-4\alpha-4})/2<2(\alpha-1)$.

In summary, the log-concavity condition $s_{k-1}\cdot s_{k+1}\leq s_{k}%
^{2},1\leq k\leq\alpha-1$, holds for $\alpha\leq5$.
\end{proof}

A graph $G$ is called \textit{perfect} if $\chi(H)=\omega(H)$ for any induced
subgraph $H$ of $G$, where $\chi(H)$ denotes the chromatic number of $H$
(Berge, \cite{Berge0}). Lov\'{a}sz proved the theorem claiming that a graph
$G$ is perfect if and only if $\left\vert V(H)\right\vert \leq\alpha
(H)\cdot\omega(H)$ for any induced subgraph $H$ of $G$ (see \cite{Lovasz}).

\begin{proposition}
\label{prop3}If $G$ is a perfect graph with $\alpha(G)=\alpha$ and
$\omega=\omega(G)$, then $s_{\left\lceil (\omega\alpha-1)/\left(
\omega+1\right)  \right\rceil }\geq...\geq s_{\alpha-1}\geq s_{\alpha}$.
\end{proposition}

\begin{proof}
Let $S$ be a stable set in $G$ of size $k\geq0$. Then $H=G-N[S]$ is an induced
subgraph of $G$\ and has $\alpha(H)\leq\alpha-k$. Therefore, by Lov\'{a}sz's
theorem,
\[
\left\vert V(H)\right\vert \leq\omega(H)\cdot\alpha(H)\leq\omega
(H)\cdot(\alpha-k)\leq\omega\cdot(\alpha-k)
\]
and, hence, $\omega_{\alpha-k}\leq\omega\cdot(\alpha-k)$. Further, according
to Lemma \ref{lem1}, we obtain that $(k+1)\cdot s_{k+1}\leq\omega\cdot
(\alpha-k)\cdot s_{k},0\leq k<\alpha$. Now, $s_{k+1}\leq s_{k}$ is true while
$k+1\geq\omega\cdot(\alpha-k)$, i.e., for $k\geq(\omega\alpha-1)/\left(
\omega+1\right)  $.
\end{proof}

In fact, in Proposition \ref{prop3} there is some $k$ such that $\left\lceil
(\omega\alpha-1)/\left(  \omega+1\right)  \right\rceil \leq k<\alpha$ if and
only if $\alpha-\frac{1+\alpha}{1+\omega}\leq\alpha-1$, i.e., for $\alpha
\geq\omega$. It is worth mentioning that, for general graphs,\emph{\ }Lemma
\ref{lem1} assures that if a graph $G$ satisfies $\omega(G)\leq\alpha
=\alpha(G)$, then $s_{\alpha}\leq s_{\alpha-1}$. However, the inverse
assertion is not true, e.g., $\alpha(K_{4}-e)=2<3=\omega(K_{4}-e)$ and
$I(K_{4}-e;x)=1+4x+x^{2}$, where by $K_{4}-e$ we mean the graph obtained from
$K_{4}$ by deleting one of its edges.

For non-perfect graphs, Proposition \ref{prop3} is not necessarily false, for
example, $I(C_{7};x)=$ $1+7x+14x^{2}+7x^{3}$. However, the graph
$G=\sqcup4C_{5}$ is not perfect, $\alpha(G)=8,\omega(G)=2$ and
\begin{align*}
I(\sqcup4C_{5};x)  & =\left(  1+5x+5x^{2}\right)  ^{4}=1+20x+170x^{2}%
+800x^{3}+2275x^{4}+\\
& +4000x^{5}+\mathbf{4250}x^{6}+2500x^{7}+625x^{8}%
\end{align*}
is log-concave, but it does not satisfies Proposition \ref{prop3}, since
$\left\lceil (\omega\alpha-1)/\left(  \omega+1\right)  \right\rceil
=\left\lceil (2\cdot8-1)/\left(  2+1\right)  \right\rceil $ $=$ $5$ and
$s_{5}=4000<4250=s_{6}$.

Any minimal imperfect graph $G$, i.e., $G=C_{2n+1},n\geq2$, or $G=$
$\overline{C_{2n+1}},n\geq2$, is claw-free and, consequently, its $I(G;x)$ is
log-concave, by Hamidoune's theorem, \cite{Hamidoune}. However, there are
non-perfect graphs, whose independence polynomials are not unimodal, e.g., the
disconnected graph $G=\left(  K_{95}+\sqcup4K_{3}\right)  \sqcup C_{5}$ has%
\begin{align*}
I(G;x)  & =\left(  1+107x+54x^{2}+108x^{3}+81x^{4}\right)  \left(
1+5x+5x^{2}\right) \\
& =1+112x+594x^{2}+\mathbf{913}x^{3}+891x^{4}+\mathbf{945}x^{5}+405x^{6}.
\end{align*}
Let $H=K_{97}+\sqcup4K_{3}$, and $G$ be the graph obtained from $H$ by adding
an edge from a vertex of $K_{97}$ to a vertex of some $C_{5}$. Then $G $ is a
connected imperfect graph whose $I(G;x)$ is not unimodal, since
\begin{align*}
I(G;x)  & =\left(  1+109x+54x^{2}+108x^{3}+81x^{4}\right)  \left(
1+4x+3x^{2}\right) \\
& +x(1+2x)\left(  1+108x+54x^{2}+108x^{3}+81x^{4}\right) \\
& =1+114x+603x^{2}+\mathbf{921}x^{3}+891x^{4}+\mathbf{945}x^{5}+405x^{6}.
\end{align*}

Since any bipartite graph $G$ is perfect and has $\omega(G)\leq2$, we obtain
the following result.

\begin{corollary}
\label{cor2}If $G$ is a bipartite graph with $\alpha(G)=\alpha\geq1$,

then $s_{\left\lceil (2\alpha-1)/3\right\rceil }\geq...\geq s_{\alpha-1}\geq
s_{\alpha}$.
\end{corollary}

In particular, we infer a similar result for trees, whose importance is
significant vis-\`{a}-vis the conjecture of Alavi \textit{et al}.

\begin{corollary}
If $T$ is a tree with $\alpha(T)=\alpha$, then $s_{\left\lceil (2\alpha
-1)/3\right\rceil }\geq...\geq s_{\alpha-1}\geq s_{\alpha}$.
\end{corollary}

For non-bipartite graphs, Corollary \ref{cor2} is not necessarily false (see
the graphs in Figure \ref{fig22}).

\section{Conclusions}

In this paper we prove that for very well-covered graphs the "chaotic
interval" $(s_{\left\lceil \alpha/2\right\rceil },s_{\left\lceil
\alpha/2\right\rceil +1},...,s_{\alpha})$ involved in the roller-coaster
conjecture of Michael and Traves can be shorten to $(s_{\left\lceil
\alpha/2\right\rceil },s_{\left\lceil \alpha/2\right\rceil +1}%
,...,s_{\left\lceil (2\alpha-1)/3\right\rceil })$. It seems that one can get
even deeper results, by using more efficiently the power of the new defined
parameter $\omega_{k}$.

We also conclude with the two following conjectures sharpening the conjectures
of Brown \textit{et al. }and Alavi \textit{et al. }respectively.

\begin{conjecture}
$I(G;x)$ is log-concave for any very well-covered graph $G$.
\end{conjecture}

\begin{conjecture}
$I(T;x)$ is log-concave for any (well-covered) tree $T$.
\end{conjecture}

\end{document}